\documentclass[12pt,a4paper]{article}
\usepackage{amsmath,amsfonts,amssymb,amsthm,amscd}
\newtheorem{thm}{Theorem}
\newtheorem{lem}[thm]{Lemma}
\newtheorem{prop}[thm]{Proposition}

\newtheorem{cor}[thm]{Corollary}

\newtheorem{notations}[thm]{Notations}

\begin{document}

\author{Liana David}

\title{The fundamental form of almost-quaternionic Hermitian manifolds}

\maketitle

{\it Author's address:}  Institute of Mathematics "Simion Stoilow"
of the Romanian Academy, Calea Grivitei nr 21, Sector 1,
Bucharest, Romania; tel. 0040-21-3196531; fax 0040-21-3196505;
e-mail
address: liana.david@imar.ro  and liana.r.david@gmail.com\\

\textbf{Abstract:} We prove that if the fundamental $4$-form
$\Omega$ of an almost-quaternionic Hermitian manifold $(M, Q, g)$
of dimension $4n\geq 8$ satisfies the conformal-Killing
equation, then $(M, Q, g)$ is quaternionic-K\"{a}hler.\\


{\it Subject Classification:} Real and complex differential geometry; Geometric theory
of differential equations.\\


\section{Introduction}

Conformal-Killing (respectively, Killing) $1$-forms are dual to
conformal-Killing (respectively, Killing) vector fields. More
generally, a $p$-form $\psi$ ($p\geq 1$) on a Riemannian manifold
$(M^{m}, g)$ is conformal-Killing, if it satisfies the
conformal-Killing equation
\begin{equation}\label{confkil}
\nabla_{X} \psi  =\frac{1}{p+1} i_{X}d\psi  -\frac{1}{m-p+1}X\land
\delta \psi ,\quad\forall X\in TM,
\end{equation}
where $\nabla$ is the Levi-Civita connection and (like everywhere
in this note) we identify tangent vectors with $1$-forms by means
of the Riemannian duality. Co-closed conformal-Killing forms are
called Killing. Note that $\psi$ is Killing if and only if its
covariant derivative is totally skew, or, equivalently,
$(\nabla_{X}\psi )(X,\cdot ) =0$ for any vector field $X$.

Conformal-Killing forms exist on spaces of constant curvature, on
Sasaki manifolds \cite{sem} and on some classes of K\"{a}hler
manifolds, like Bochner-flat K\"{a}hler manifolds and
conformally-Einstein K\"{a}hler manifolds \cite{ap}, \cite{mor}.
On compact quaternionic-K\"{a}hler manifolds of dimension at least
eight, there are no non-parallel conformal-Killing $2$-forms,
unless the quaternionic-K\"{a}hler manifold is isomorphic to the
standard quaternionic projective space, in which case the space of
conformal-Killing $2$-forms is naturally isomorphic to the space
of Killing vector fields \cite{david}.

Conformal-Killing forms exist also on manifolds which admit
twistor spinors \cite{sem}. Recall that a twistor spinor on a
Riemannian spin manifold $(M^{m}, g)$ is a section $\rho$ of the
spinor bundle, which satisfies the equation $\nabla_{X}\rho =
-\frac{1}{m} X\cdot D\rho$, where $X$ is any vector field, $D$ is
the Dirac operator and "$\cdot $" denotes the Clifford
multiplication. If $\rho_{1}$ and $\rho_{2}$ are twistor spinors,
then the $p$-form
$$
\omega_{p}(X_{1},\cdots ,  X_{p}) =\langle  (X_{1}\wedge \cdots
\wedge X_{p})\cdot \rho_{1}, \rho_{2}\rangle
$$
is conformal-Killing (for any $p\geq 1$). For a survey on
conformal-Killing forms, see for example
\cite{sem}.\\

The starting point of this note is a result proved in \cite{sem},
which states that if the K\"{a}hler form of an almost-Hermitian
manifold is conformal-Killing, then the almost-Hermitian manifold
is nearly K\"{a}hler. Our main Theorem is an analogue of this
result in quaternionic geometry and is stated as follows:

\begin{thm}\label{qkmain} Let $(M^{4n},Q, g)$ be an almost-quaternionic Hermitian
manifold, of dimension $4n\geq 8.$ Suppose that the fundamental
$4$-form $\Omega$ of $(M, Q,g)$ is conformal-Killing. Then $(M, Q,
g)$ is quaternionic-K\"{a}hler.
\end{thm}

Theorem \ref{qkmain} generalizes a result proved in \cite{swann3},
namely that in dimension at least eight, a nearly
quaternionic-K\"{a}hler manifold (i.e. an almost-quaternionic
Hermitian manifold for which the fundamental $4$-form is a Killing
form) is necessarily quaternionic-K\"{a}hler.

The paper is organized as follows: in Section \ref{section1} we
recall basic facts on quaternionic Hermitian geometry. Section
\ref{section2} is devoted to the proof of our main result, which
is based on a representation theoretic argument. Similar arguments
were already employed in
\cite{swann2} and \cite{swann3}.\\

\section{Quaternionic Hermitian geometry}\label{section1}

Let $M$ be a manifold of dimension $4n\geq 8$ (in all our
considerations the dimension of the manifold will be at least
eight). An almost-quaternionic structure on $M$ is a rank-three
vector sub-bundle $Q\subset \mathrm{End}(TM)$, locally generated
by three anti-commuting almost complex structures $\{ J_{1},
J_{2}, J_{3}\}$ which satisfy $J_{1}\circ J_{2} = J_{3}.$ Such a
triple of almost complex structures is usually called a (local)
admissible basis of $Q$. An almost-quaternionic Hermitian
structure on $M$ consists of an almost-quaternionic structure $Q$
and a Riemannian metric $g$ compatible with $Q$, which means that
$$
g(JX, JY) = g (X, Y),\quad\forall J\in Q, \quad
J^{2}=-\mathrm{Id}, \quad\forall X, Y\in TM.
$$
In the language of $G$-structures, an almost-quaternionic
Hermitian structure on a $4n$-dimensional manifold is an
$Sp(n)Sp(1)$-structure. Therefore, on an almost-quaternionic
Hermitian manifold $(M^{4n}, g, Q)$ there are two locally defined
complex vector bundles $E$ and $H$, of rank $2n$ and $2$
respectively, associated to the standard representations of
$Sp(n)$ and $Sp(1)$ on $\mathbb{E}=\mathbb{C}^{2n}$ and
$\mathbb{H}=\mathbb{C}^{2}$. Let $\omega_{E}\in\Lambda^{2}(E^{*})$
and $j_{E}:E\rightarrow E$ be the standard symplectic form and
quaternionic structure of the bundle $E$, defined by the
$Sp(n)$-invariant complex symplectic form and quaternionic
structure of $\mathbb{E}.$ We shall often identify $E$ with
$E^{*}$ by means of the map $e\rightarrow \omega_{E}(e,\cdot )$,
so that $\omega_{E}$ will sometimes be considered as a bivector on
$E$. For any $r\geq 2$ we shall denote by $\Lambda^{r}_{0}E\subset
\Lambda^{r}E$ the kernel of the natural contraction
\begin{equation}\label{omegacon}
\omega_{E}\bullet :\Lambda^{r}E\rightarrow \Lambda^{r-2}E
\end{equation}
with the symplectic form $\omega_{E}$, defined by
$$
\omega_{E}\bullet (e_{1}\wedge\cdots \wedge e_{r}) =
\sum_{i<j}(-1)^{i+j+1}\omega_{E}(e_{i}, e_{j}) e_{1}\wedge\cdots
\wedge \widehat{{e}_{i}}\wedge\cdots \wedge
\widehat{{e}_{j}}\wedge\cdots \wedge e_{r}
$$
where the hat denotes that the term is omitted. By means of
contraction and wedge product with $\omega_{E}$ we can decompose
$\Lambda^{r}E$ as
\begin{equation}\label{30}
\Lambda^{r}E = \Lambda^{r}_{0}E\oplus \omega_{E}\wedge
\Lambda^{r-2}_{0}E\oplus
\omega_{E}^{2}\wedge\Lambda^{r-4}_{0}E\oplus\cdots
\end{equation}
The map $j_{E}$ is complex anti-linear and
$$
j_{E}^{2}= -\mathrm{Id},\quad\omega_{E}(j_{E}u, j_{E}v) =
\overline{\omega_{E}(u, v)},\quad \omega_{E}(e, j_{E}{e})>0,
$$
for any $u, v\in E$ and $e\in E\setminus \{ 0\}$. To simplify
notations, for a vector $e\in E$ we shall often denote $\tilde{e}
:= j_{E}(e)$ its image through the quaternionic structure of $E$.
Similar conventions will be used for the standard symplectic form
$\omega_{H}\in\Lambda^{2}(H^{*})$ and quaternionic structure
$j_{H}:H\rightarrow H$ of the
bundle $H.$\\

The bundles $E$ and $H$ play the role of spin bundles from
conformal geometry. In particular,
\begin{equation}\label{tang}
T_{\mathbb{C}}M = E\otimes_{\mathbb{C}} H
\end{equation}
and the complex bilinear extension of the Riemannian metric $g$ to
$T_{\mathbb{C}}M$ is the tensor product
$\omega_{E}\otimes\omega_{H}.$ Decomposition (\ref{tang}) induces
decompositions of the form bundles in any degree. In particular,
the bundles of $2$ and $3$-forms decompose as (see \cite{sal})
\begin{equation}\label{bivectors}
\Lambda^{2}(T_{\mathbb{C}}M) = S^{2}H\oplus S^{2}E\oplus
S^{2}H\Lambda^{2}_{0}E
\end{equation}
\begin{equation}\label{lista1}
\Lambda^{3}(T_{\mathbb{C}}M)= H(E\oplus K)\oplus S^{3}H
(\Lambda^{3}_{0}E\oplus E).
\end{equation}
(In (\ref{bivectors}) and (\ref{lista1}), and often in this note,
we omit the tensor product signs). In (\ref{bivectors}) $S^{2}H$
and $S^{2}E$ are complexifications of the bundle $Q$ and,
respectively, of the bundle of $Q$-Hermitian $2$-forms, i.e.
$2$-forms $\psi\in\Lambda^{2}(T^{*}M)$ which satisfy
$$
\psi (JX, JY) = \psi (X, Y),\quad\forall J\in Q,\quad J^{2} =
-\mathrm{Id}, \quad\forall X, Y\in TM.
$$
In (\ref{lista1}) $K$ denotes the vector bundle associated to the
$Sp(n)$-module $\mathbb{K}$, which arises into the irreducible
decomposition
\begin{equation}\label{k}
\mathbb{E}\otimes
\Lambda^{2}_{0}\mathbb{E}\cong\Lambda^{3}_{0}\mathbb{E} \oplus
\mathbb{E}\oplus \mathbb{K}
\end{equation}
under the action of $Sp(n).$ A vector from $\mathbb{E}\otimes
\Lambda^{2}_{0}\mathbb{E}$ has non-trivial component on
$\mathbb{K}$ if and only if it is not totally skew.

\begin{notations}{\rm We shall identify bundles with their
complexification, without additional explanations. For example, in
(\ref{bivectors}) $S^{2}H\Lambda^{2}_{0}E$ is a complex sub-bundle
of $\Lambda^{2}(T_{\mathbb{C}}M)$. We shall use the same notation
for its real part, which is a sub-bundle of $\Lambda^{2}(TM).$}
\end{notations}

An almost-quaternionic Hermitian manifold $(M,g, Q)$ has a
canonical $4$-form, defined, in terms of an arbitrary admissible
basis $\{ J_{1}, J_{2}, J_{3}\}$ of $Q$, by
$$
\Omega = \omega_{1}\wedge\omega_{1}+\omega_{2}\wedge\omega_{2}
+\omega_{3}\wedge\omega_{3},
$$
where $\omega_{i}:= g(J_{i}\cdot , \cdot )$ are the K\"{a}hler
forms corresponding to $(g, J_{i}).$ As proved in \cite{bryant}
and \cite{swann2}, the covariant derivative $\nabla \Omega$ with
respect to the Levi-Civita connection $\nabla$ of $g$ is a section
of $T^{*}M\otimes (S^{2}H\Lambda^{2}_{0}E)$, where
$S^{2}H\Lambda^{2}_{0}E$ is embedded into $\Lambda^{4}(T^{*}M)$
(identified with $\Lambda^{4}(TM)$ using the Riemannian metric),
in the following way. Note first that $\Lambda^{2}(S^{2}H)$ is
canonically isomorphic to $S^{2}H$ (this is because $S^{2}H$ is
the complexification of $Q$, which has a natural metric and
orientation, for which any admissible basis $\{ J_{1}, J_{2},
J_{3}\}$ is orthonormal and positively oriented). The map
\begin{equation}\label{explicatii}
S^{2}H\Lambda^{2}_{0}E\cong
\Lambda^{2}(S^{2}H)\Lambda^{2}_{0}E\rightarrow
\Lambda^{4}_{\mathbb{C}}(TM)
\end{equation}
defined by
\begin{equation}\label{s12}
(s_{1}\wedge s_{2})\beta\rightarrow s_{1}\beta\wedge
s_{2}\omega_{E} - s_{2}\beta\wedge s_{1}\omega_{E}, \quad \forall
s_{1}, s_{2}\in S^{2}H,\quad\forall \beta\in \Lambda^{2}_{0}E
\end{equation}
is the promised embedding of $S^{2}H\Lambda^{2}_{0}E$ into
$\Lambda^{4}(TM)$.

An almost-quaternionic Hermitian manifold $(M, Q,g)$ is
quaternionic-K\"{a}hler if the Levi-Civita connection $\nabla$ of
$g$ preserves the bundle $Q$, or, equivalently, the fundamental
$4$-form $\Omega$ is parallel with respect to $\nabla$. In fact,
as already mentioned in the Introduction, according to Theorem 1.2
of \cite{swann3} the weaker condition $(\nabla_{X}\Omega )(X,\cdot
)=0$, for any vector field $X$, implies that $(M, Q, g)$ is
quaternionic-K\"{a}hler.

\section{Proof of the main result}\label{section2}

In this Section we prove our main result. Let $(M, Q, g)$ be an
almost-quaternionic Hermitian manifold, whose fundamental $4$-form
$\Omega$ is conformal-Killing. In order to prove that $\Omega$ is
parallel with respect to the Levi-Civita connection $\nabla$, it
is enough to show that it is co-closed (being conformal-Killing,
$\Omega$ is co-closed if and only if it is Killing, if and only if
it is parallel, by Theorem 1.2 of \cite{swann3} already mentioned
before). Recall now that $\nabla\Omega$ is a section of
$T^{*}M\otimes (S^{2}H\Lambda^{2}_{0}E)$, which decomposes into
irreducible sub-bundles as
\begin{equation}\label{lista}
T^{*}_{\mathbb{C}}M\otimes (S^{2}H\Lambda^{2}_{0}E)= HE\oplus
H\Lambda^{3}_{0}E\oplus HK\oplus (S^{3}H)E\oplus
S^{3}H\Lambda^{3}_{0}E\oplus (S^{3}H)K.
\end{equation}
Decomposition (\ref{lista}) follows from (\ref{k}), together with
the irreducible decomposition
$$
\mathbb{H}\otimes S^{2}\mathbb{H}\cong S^{3}\mathbb{H}\oplus
\mathbb{H}
$$
of $\mathbb{H}\otimes S^{2}\mathbb{H}$ under $Sp(1).$ While
$H\Lambda^{3}_{0}E$ and $(S^{3}H)K$ are irreducible sub-bundles of
$T^{*}_{\mathbb{C}}M\otimes (S^{2}H\Lambda^{2}_{0}E)$, see
(\ref{lista}), they are not irreducible sub-bundles of
$\Lambda^{3}(T_{\mathbb{C}}M)$, see (\ref{lista1}). These
observations readily imply that if $\nabla\Omega$ is a section of
$H\Lambda^{3}_{0}E\oplus (S^{3}H)K$, then $\Omega$ is co-closed:
just write $\delta\Omega = - \sum_{i} (\nabla_{E_{i}}\Omega
)(E_{i}, \cdot )$, where $\{ E_{i}\}$ is a local orthonormal frame
of $TM$, and use the fact that an invariant linear map between
non-isomorphic irreducible representations is identically zero.
(Actually, by Theorem 2.3 of \cite{swann3}, also the converse is
true: if $\delta \Omega =0$ then $\nabla\Omega$ is a section of
$H\Lambda^{3}_{0}E\oplus (S^{3}H)K$).

Therefore, we aim to show that $\nabla\Omega$ is a section of
$H\Lambda^{3}_{0}E\oplus (S^{3}H)K$. For this, we define the
algebraic conformal-Killing operator
$$
\mathcal T : T^{*}M\otimes \Lambda^{4}(TM)\rightarrow
T^{*}M\otimes \Lambda^{4}(TM),
$$
by
\begin{equation}\label{detf}
\mathcal T (\gamma \otimes\alpha )(X) = \frac{4}{5}\gamma
(X)\alpha +\frac{1}{5}\gamma\wedge i_{X}\alpha   -
\frac{1}{4n-3}X\wedge i_{\gamma}\alpha
\end{equation}
where $\gamma\in T^{*}M$ (is identified with a vector using the
Riemannian metric), $\alpha\in\Lambda^{4}(TM)$ and $X\in TM$. Note
that, for any $4$-form $\psi\in \Omega^{4}(M)$,
\begin{equation}\label{adaugat}
{\mathcal T}(\nabla \psi )(X) = \nabla_{X}\psi
-\frac{1}{5}i_{X}d\psi +\frac{1}{4n-3}X\wedge \delta \psi
,\quad\forall X\in TM.
\end{equation}
In particular, since $\Omega$ is conformal-Killing,
\begin{equation}\label{kernel}
{\mathcal T}(\nabla\Omega )=0.
\end{equation}
The operator $\mathcal T$ is $Sp(n)Sp(1)$-invariant and we extend
it, by complex linearity, to $T^{*}_{\mathbb{C}}M\otimes
\Lambda^{4}(T_{\mathbb{C}}M)$. Define
$$
{\mathcal S}:= T^{*}_{\mathbb{C}}M\otimes
(S^{2}H\Lambda^{2}_{0}E)\ominus \left( H\Lambda^{3}_{0}E\oplus
(S^{3}H)K\right).
$$
From (\ref{lista}), the irreducible sub-bundles of $\mathcal S$
are
\begin{equation}\label{ll}
HE,\quad HK ,\quad (S^{3}H)E,\quad S^{3}H\Lambda^{3}_{0}E.
\end{equation}
For any irreducible sub-bundle $W$ of $\mathcal S$, we will
determine an $Sp(n)Sp(1)$-invariant linear map
$$
{\mathcal T}_{W}:T^{*}_{\mathbb{C}}M\otimes
\Lambda^{4}(T_{\mathbb{C}}M) \rightarrow W
$$
which factors through ${\mathcal T}$ (i.e. ${\mathcal T}_{W}
=\mathrm{pr}_{W}\circ{\mathcal T}$ is the composition of
${\mathcal T}$ with an $Sp(n)Sp(1)$-invariant linear map
$\mathrm{pr}_{W}$ from $T^{*}_{\mathbb{C}}M\otimes
\Lambda^{4}(T_{\mathbb{C}}M)$ to $W$) such that the restriction of
${\mathcal T}_{W}$ to $T^{*}_{\mathbb{C}}M\otimes
(S^{2}H\Lambda^{2}_{0}E)$ is non-zero. An easy argument which uses
(\ref{kernel}),  Schur's Lemma and the fact that irreducible
sub-bundles of $T^{*}_{\mathbb{C}}M\otimes
(S^{2}H\Lambda^{2}_{0}E)$ are pairwise non-isomorphic, would then
imply that $\nabla\Omega$ has trivial component on $W$ and
therefore that $\nabla\Omega$ is a section
of $H\Lambda^{3}_{0}E\oplus (S^{3}H)K$, as needed.\\

In order to define the maps ${\mathcal T}_{W}$, we apply several
suitable contractions to the algebraic conformal-Killing operator
${\mathcal T}$. We first define ${\mathcal T}_{HE}$ and ${\mathcal
T}_{HK}$ as follows. For a section $\eta$ of
$T^{*}_{\mathbb{C}}M\otimes\Lambda^{4}(T_{\mathbb{C}}M)$, define
$\omega_{E}\bullet {\mathcal T}(\eta )$, a $1$-form with values in
$(S^{2}H)\Lambda^{2}(T_{\mathbb{C}}M)$, by
\begin{equation}\label{oe}
\omega_{E}\bullet \left( {\mathcal T}(\eta )\right)
(X):=\omega_{E}\bullet \left( {\mathcal T}(\eta )(X)\right)
,\quad\forall X\in TM,
\end{equation}
where in (\ref{oe}) ${\mathcal T}(\eta )(X)$ belongs to
$\Lambda^{4}(T_{\mathbb{C}}M)$ (is the value of the
$\Lambda^{4}(T_{\mathbb{C}}M)$-valued $1$-form ${\mathcal T}(\eta
)$ on $X\in T_{\mathbb{C}}M$) and
\begin{equation}\label{contraction}
\omega_{E}\bullet : \Lambda^{4}(T_{\mathbb{C}}M)\rightarrow
(S^{2}H)\Lambda^{2}(T_{\mathbb{C}}M)
\end{equation}
denotes the contraction with $\omega_{E}$, which on decomposable
multi-vectors
$$
\beta = h_{1}e_{1}\wedge\cdots \wedge h_{4}e_{4}\in
\Lambda^{4}(T_{\mathbb{C}}M)
$$
takes value
$$
\omega_{E}(\beta )=
\sum_{i<j}(-1)^{i+j+1}\omega_{E}(e_{i},e_{j})(h_{i}h_{j}+h_{j}h_{i})
h_{1}e_{1}\wedge\cdots \wedge \widehat{h_{i}e_{i}}\wedge\cdots
\wedge \widehat{h_{j}e_{j}}\wedge\cdots \wedge h_{4}e_{4}.
$$
Next, we define $\omega_{H}\bullet\omega_{E}\bullet {\mathcal
T}(\eta )$, by contracting $\omega_{E}\bullet {\mathcal T}(\eta
)$, which is a section of $HE\otimes
(S^{2}H)\Lambda^{2}(T_{\mathbb{C}}M)$, with $\omega_{H}$ in the
first two $H$-variables. Therefore,
$\omega_{H}\bullet\omega_{E}\bullet {\mathcal T}(\eta )$ is a
section of $EH \Lambda^{2}(T_{\mathbb{C}}M)$. Considering $EH
\Lambda^{2}(T_{\mathbb{C}}M)$ naturally embedded into $EH(HHEE)$,
we contract further $\omega_{H}\bullet\omega_{E}\bullet {\mathcal
T}(\eta )$ with $\omega_{H}$ again in the first two $H$-variables.
The result is a section $\omega_{H}^{2}\bullet\omega_{E}\bullet
{\mathcal T}(\eta )$ of $HE E E$. Applying suitable projections to
$\omega_{H}^{2}\bullet\omega_{E}\bullet {\mathcal T}(\eta )$ we
finally obtain ${\mathcal T}_{HE}(\eta )$ and ${\mathcal
T}_{HK}(\eta ) $, as follows.

The contraction of $\omega_{H}^{2}\bullet\omega_{E}\bullet
{\mathcal T}(\eta )$ with $\omega_{E}$ in the first two
$E$-variables defines
\begin{equation}\label{01}
{\mathcal T}_{HE}(\eta ):= \omega_{E}\bullet \omega_{H}^{2}\bullet
\omega_{E}\bullet {\mathcal T}(\eta ) .
\end{equation}
Similarly, we can project $\omega_{H}^{2}\bullet\omega_{E}\bullet
{\mathcal T}(\eta )$ to $H\otimes E\Lambda^{2}_{0}E$ and then to
$HK$, by means of the decomposition (\ref{k}) (translated to
vector bundles). The result of this projection is the value of
${\mathcal T}_{HK}$ on $\eta$. More precisely,
\begin{equation}\label{02}
{\mathcal T}_{HK}(\eta ):=
\mathrm{pr}_{HK}\left(\omega_{H}^{2}\bullet \omega_{E}\bullet
{\mathcal T}(\eta )\right) .
\end{equation}

\begin{prop}\label{prop1} The operators ${\mathcal T}_{HE}$ and
${\mathcal T}_{HK}$ defined by (\ref{01}) and (\ref{02}) are
non-trivial on $T^{*}_{\mathbb{C}}M\otimes
(S^{2}H\Lambda^{2}_{0}E)$.
\end{prop}

In order to prove Proposition \ref{prop1}, we will show that
${\mathcal T}_{HE}$ and ${\mathcal T}_{HK}$ take non-zero value on
$\gamma_{0}\alpha_{0}$, where
\begin{equation}\label{alphag}
\gamma_{0} := \tilde{e}_{1}h,\quad \alpha_{0} := e_{1}h\wedge
e_{2}h\wedge e_{i}\tilde{h}\wedge \tilde{e}_{i}\tilde{h} -
e_{1}\tilde{h}\wedge e_{2}\tilde{h}\wedge e_{i}{h}\wedge
\tilde{e}_{i}h
\end{equation}
was already considered in \cite{swann3}. In (\ref{alphag}) $\{
e_{1}, \cdots , e_{2n}\}$ is a unitary basis of (local) sections
of $E$, with respect to the (positive definite) Hermitian metric
$g_{E}:= \omega_{E}(\cdot , j_{E}\cdot )$, chosen such that
$e_{n+j}= \tilde{e}_{j}$ for any $1\leq j\leq n$, and $\{ h,
\tilde{h}\} $ is a unitary basis of (local) sections of $H$, with
respect to $g_{H}:=\omega_{H}(\cdot , j_{H}\cdot )$. In order to
simplify notations, in (\ref{alphag}) and bellow we omit the
summation sign over $1\leq i\leq 2n$. The symplectic forms of $E$
and $H$ can be written as
\begin{equation}\label{ba}
\omega_{E}=\frac{1}{2}e_{i}\wedge\tilde{e}_{i}\in\Lambda^{2}E,
\quad \omega_{H} = h\wedge\tilde{h}\in\Lambda^{2}H.
\end{equation}
From (\ref{s12}) and (\ref{ba}), $\alpha_{0}$ is a section of the
sub-bundle $S^{2}H\Lambda^{2}_{0}E$ of
$\Lambda^{4}(T_{\mathbb{C}}M)$ and $\gamma_{0}\alpha_{0}$ is a
section of $T^{*}_{\mathbb{C}}M\otimes
(S^{2}H\Lambda^{2}_{0}E).$\\

We divide the proof of Proposition \ref{prop1} into the following
two Lemmas.

\begin{lem}\label{lema1} The section
$\mathrm{pr}_{HE\Lambda^{2}_{0}E} \left(\omega_{H}^{2}\bullet
\omega_{E}\bullet {\mathcal T}(\gamma_{0}\alpha_{0})\right)$ is
not totally skew in the $E$-variables. In particular, ${\mathcal
T}_{HK} (\gamma_{0}\alpha_{0})\neq 0$.
\end{lem}

\begin{proof} A straightforward computation shows that
$$
i_{\gamma_{0}}\alpha_{0}= e_{i}h\wedge\tilde{e}_{i}h\wedge
e_{2}\tilde{h} -2 e_{1}h\wedge e_{2}h\wedge\tilde{e}_{1}\tilde{h}.
$$
Therefore, using (\ref{detf}), we can write
\begin{equation}\label{t}
{\mathcal T}(\gamma_{0}\alpha_{0}) =
\frac{4}{5}\gamma_{0}\alpha_{0} +\frac{1}{5}\gamma_{0}\wedge
\alpha_{0} (\cdot ) -\frac {1}{4n-3} (F-2G),
\end{equation}
where $\gamma_{0}\wedge\alpha_{0} (\cdot )$ is a $1$-form with
values in $\Lambda^{4}(T_{\mathbb{C}}M)$, whose natural
contraction with a vector $X\in T_{\mathbb{C}}M$ is
$\gamma_{0}\wedge i_{X}\alpha_{0}.$ Similarly, $F$ and $G$ are
defined by
\begin{align*}
F(X)&:= X\wedge e_{i}h\wedge \tilde{e}_{i}h\wedge e_{2}\tilde{h}\\
G(X)&:= X\wedge e_{1}h\wedge e_{2}h\wedge \tilde{e}_{1}\tilde{h}.
\end{align*}
Now, it is straightforward to check  that
\begin{align*}
\omega_{H}^{2}\bullet\omega_{E}\bullet (\gamma_{0}\alpha_{0}
) &= -4n h( \tilde{e}_{1}e_{1}e_{2} - \tilde{e}_{1}e_{2}e_{1})\\
\omega_{H}^{2}\bullet\omega_{E}\bullet \left( \gamma_{0}\wedge
\alpha_{0} (\cdot )\right) &= 2h( -e_{i}\tilde{e}_{i}e_{2}
+\tilde{e}_{1}e_{1}e_{2} -\tilde{e}_{1}e_{2}e_{1} +\tilde{e}_{i}e_{i}e_{2})\\
& + h(-e_{2}e_{i}\tilde{e}_{i}+e_{2}\tilde{e}_{i}e_{i}
+e_{i}e_{2}\tilde{e}_{i} -\tilde{e}_{i} e_{2}e_{i})\\
&+(4n+2)h( e_{2}\tilde{e}_{1}e_{1} - e_{1}\tilde{e}_{1}e_{2})\\
&+4h (e_{1}e_{2}\tilde{e}_{1}-e_{2}e_{1}\tilde{e}_{1})
\end{align*}
and also
\begin{align*}
\omega_{H}^{2}\bullet\omega_{E}\bullet F&= -(4n-4)h
e_{i}e_{2}\tilde{e}_{i} + 3h
(e_{2}e_{i}\tilde{e}_{i}-e_{2}\tilde{e}_{i}e_{i})\\
\omega_{H}^{2}\bullet\omega_{E}\bullet G&= 3h
(e_{1}\tilde{e}_{1}e_{2} - e_{2}\tilde{e}_{1}e_{1}
-\tilde{e}_{1}e_{2}e_{1} +\tilde{e}_{1}e_{1}e_{2}) - h
e_{i}e_{2}\tilde{e}_{i}.\\
\end{align*}
These relations combined with (\ref{t}) readily imply that
\begin{align*}
\omega_{H}^{2}\bullet \omega_{E}\bullet {\mathcal
T}(\gamma_{0}\alpha_{0})&= \lambda_{1}h\tilde{e}_{1}( e_{1}\wedge
e_{2})+\lambda_{2}h(e_{2}\tilde{e}_{1}e_{1}-e_{1}\tilde{e}_{1}e_{2})\\
&+\lambda_{3}he_{2}(
e_{i}\wedge\tilde{e}_{i})+\lambda_{4}he_{i}e_{2}\tilde{e}_{i}\\
&+\frac{h}{5}\left( 4(e_{1}\wedge e_{2}) \tilde{e}_{1}
+2(\tilde{e}_{i}\wedge e_{i}) e_{2}
-\tilde{e}_{i}e_{2}e_{i}\right) ,
\end{align*}
with constants
\begin{align*}
\lambda_{1}=\frac{8(-8n^{2}+7n+3)}{5(4n-3)},\quad \lambda_{2}=
\frac{4(4n^{2}-n-9)}{5(4n-3)},\quad
\lambda_{3}=-\frac{4(n+3)}{5(4n-3)}
\end{align*}
and
$$
\lambda_{4}=\frac{24n-33}{5(4n-3)}.
$$
Projecting the expression for $\omega_{H}^{2}\bullet
\omega_{E}\bullet {\mathcal T}(\gamma_{0}\alpha_{0})$ obtained
above onto $HE\Lambda^{2}_{0}E$ we get
\begin{align*}
\mathrm{pr}_{HE\Lambda^{2}_{0}E}\left(\omega_{H}^{2}\bullet\omega_{E}\bullet
{\mathcal T}(\gamma_{0}\alpha_{0})\right)&=
2\lambda_{1}h\tilde{e}_{1}( e_{1}\wedge
e_{2})+\left(\lambda_{2}+\frac{4}{5}\right) he_{2}(
\tilde{e}_{1}\wedge e_{1})\\
& -\left(\lambda_{2}+\frac{4}{5}\right) he_{1}(
\tilde{e}_{1}\wedge e_{2}) +\left(\lambda_{4}+\frac{2}{5}\right)
he_{i}(e_{2}\wedge \tilde{e}_{i})\\
& +\frac{3}{5} h\tilde{e}_{i}(e_{i}\wedge e_{2}) +
\frac{1}{2n}\left( \lambda_{2}-\lambda_{4} -\frac{1}{5}\right)
he_{2}(e_{i}\wedge \tilde{e}_{i}),
\end{align*}
which is not totally skew in the $E$-variables. Our claim follows.
\end{proof}

\begin{lem}\label{lema2} The value of ${\mathcal T}_{HE}$ on
$\gamma_{0}\alpha_{0}$ is
\begin{equation}\label{omegaE}
{\mathcal T}_{HE}(\gamma_{0}\alpha_{0} )  =
\frac{8n(2n+1)}{5(4n-3)}he_{2}.
\end{equation}
In particular, ${\mathcal T}_{HE}(\gamma_{0}\alpha_{0})$ is
non-zero.
\end{lem}

\begin{proof} The claim follows from a straightforward
calculation, using the expression of $\omega_{H}^{2}\bullet
\omega_{E}\bullet {\mathcal T}(\gamma_{0}\alpha_{0})$ determined
in the proof of Lemma \ref{lema1} and the definition of the
operator ${\mathcal T}_{HE}$.
\end{proof}

Lemma \ref{lema1} and Lemma \ref{lema2} conclude the proof of
Proposition \ref{prop1}.\\

We now define the maps ${\mathcal T}_{(S^{3}H)E}$ and ${\mathcal
T}_{S^{3}H\Lambda^{3}_{0}E}$. For a section $\eta$ of
$T^{*}_{\mathbb{C}}M\otimes \Lambda^{4}(T_{\mathbb{C}}M)$,
${\mathcal T}(\eta )$ is a section of
$EH\otimes\Lambda^{4}(T_{\mathbb{C}}M)$. We consider
$\omega_{H}\bullet {\mathcal T}(\eta )$, the contraction of
${\mathcal T}(\eta )$ with $\omega_{H}$ in the first two
$H$-variables, which is a section of
$EE\otimes\Lambda^{3}(T_{\mathbb{C}}M)$. Its total symmetrization
$\mathrm{sym}^{H}\left( \omega_{H}\bullet {\mathcal T}(\eta
)\right)$ in the $H$-variables is a section of
$EE(S^{3}H)\Lambda^{3}E.$ Leaving the first two $E$-variables of
$\mathrm{sym}^{H}\left( \omega_{H}\bullet {\mathcal T}(\eta
)\right)$ unchanged and contracting $\mathrm{sym}^{H}\left(
\omega_{H}\bullet {\mathcal T}(\eta )\right)$ with $\omega_{E}$ on
$\Lambda^{3}E$, as in (\ref{omegacon}), we get a section
$\omega_{E}\bullet \mathrm{sym}^{H}\left( \omega_{H}\bullet
{\mathcal T}(\eta )\right)$ of $EE(S^{3}H)E.$

To define ${\mathcal T}_{(S^{3}H)\Lambda^{3}_{0}E}(\eta ) $ and
${\mathcal T}_{(S^{3}H)E}(\eta ) $ we project $\omega_{E}\bullet
\mathrm{sym}^{H}\left( \omega_{H}\bullet {\mathcal T}(\eta
)\right)$ on $(S^{3}H)\Lambda^{3}E$ and then we project the result
on $(S^{3}H)\Lambda^{3}_{0}E$ and $(S^{3}H)E$ respectively, using
the decomposition (\ref{30}), with $r=3.$ Therefore,
\begin{equation}\label{d1}
{\mathcal T}_{(S^{3}H)\Lambda^{3}_{0}E}(\eta ):=
\mathrm{pr}_{(S^{3}H)\Lambda^{3}_{0}E} \left(\omega_{E}\bullet
\mathrm{sym}^{H}\left( \omega_{H}\bullet {\mathcal T}(\eta
)\right) \right) .
\end{equation}
Similarly,
\begin{equation}\label{d2}
{\mathcal T}_{(S^{3}H)E}(\eta ) := \omega_{E}\bullet
\mathrm{pr}_{(S^{3}H)\Lambda^{3}E}\left(\omega_{E}\bullet
\mathrm{sym}^{H}\left( \omega_{H}\bullet {\mathcal T}(\eta
)\right)\right)
\end{equation}
is the contraction of
$\mathrm{pr}_{(S^{3}H)\Lambda^{3}E}\left(\omega_{E}\bullet
\mathrm{sym}^{H}\left( \omega_{H}\bullet {\mathcal T}(\eta
)\right)\right)$ with the symplectic form $\omega_{E}.$

\begin{prop}\label{prop2} The operators
${\mathcal T}_{(S^{3}H)\Lambda^{3}_{0}E}$ and ${\mathcal
T}_{(S^{3}H)E}$ defined by (\ref{d1}) and (\ref{d2})  are
non-trivial on $T^{*}_{\mathbb{C}}M\otimes \left(
S^{2}H\Lambda^{2}_{0}E\right)$.
\end{prop}

Like in the proof of Proposition \ref{prop1}, we will show that
${\mathcal T}_{(S^{3}H)\Lambda^{3}_{0}E}(\gamma_{0}\alpha_{0})$
and ${\mathcal T}_{(S^{3}H)E}(\gamma_{0}\alpha_{0})$ are non-zero.
This is a consequence of the next Lemma.

\begin{lem}\label{lema3} The following fact holds:
\begin{align*}
\mathrm{pr}_{S^{3}H\Lambda^{3}E}\left(\omega_{E}\bullet
\mathrm{sym}^{H}\left( \omega_{H}\bullet{\mathcal
T}(\gamma_{0}\alpha_{0} )\right) \right) = -\frac{6(n-1)}{4n-3}
\mathrm{sym}^{H}(hh\tilde{h})( e_{i}\wedge\tilde{e}_{i}\wedge
e_{2})\\
-\frac{4(4n^{2}-3n+3)}{4n-3}\mathrm{sym}^{H}(hh\tilde{h})
(e_{1}\wedge e_{2}\wedge  \tilde{e}_{1}).\\
\end{align*}
\end{lem}

\begin{proof} The proof goes as in Lemma \ref{lema1}.
Applying definitions, we get:
\begin{align*}
\omega_{E}\bullet\mathrm{sym}^{H} \left(\omega_{H}\bullet
(\gamma_{0}\alpha_{0} )\right)&= 2n\mathrm{sym}^{H}(hh\tilde{h})
(\tilde{e}_{1}e_{2}e_{1}- \tilde{e}_{1}e_{1}e_{2} )\\
\omega_{E}\bullet\mathrm{sym}^{H} \left(\omega_{H}\bullet
(\gamma_{0}\wedge\alpha_{0}(\cdot ))\right) &=
(2n-4)\mathrm{sym}^{H}(hh\tilde{h})( e_{1}\tilde{e}_{1}e_{2} - e_{2}\tilde{e}_{1}e_{1})\\
&+4\mathrm{sym}^{H}(hh\tilde{h})( \tilde{e}_{1}e_{2}e_{1}-\tilde{e}_{1}e_{1}e_{2} )\\
&-2\mathrm{sym}^{H}(hh\tilde{h})(\tilde{e}_{i}e_{2}e_{i}+e_{2}\tilde{e}_{i}e_{i})\\
&+2\mathrm{sym}^{H}(hh\tilde{h})(e_{2}e_{i}\tilde{e}_{i}
+e_{i}e_{2}\tilde{e}_{i})\\
\omega_{E}\bullet\mathrm{sym}^{H} \left(\omega_{H}\bullet
F\right)& = \mathrm{sym}^{H}(hh\tilde{h})\left(
(4n-5)e_{i}\tilde{e}_{i}e_{2}
-(2n-3)e_{i}e_{2}\tilde{e}_{i}\right)\\
& +\mathrm{sym}^{H}(hh\tilde{h})\left( \tilde{e}_{i}e_{i}e_{2}-
e_{2}e_{i}\tilde{e}_{i} + e_{2}\tilde{e}_{i}e_{i}
-\tilde{e}_{i}e_{2}e_{i}\right)\\
\omega_{E}\bullet\mathrm{sym}^{H} \left(\omega_{H}\bullet
G\right)& = \mathrm{sym}^{H}(hh\tilde{h})\left(
-2e_{i}\tilde{e}_{i}e_{2} +e_{2}e_{1}\tilde{e}_{1}
-\tilde{e}_{1}e_{1}e_{2} -e_{1}e_{2}\tilde{e}_{1}\right)\\
& +\mathrm{sym}^{H}(hh\tilde{h}) \left(\tilde{e}_{1}e_{2}e_{1}+
e_{i}e_{2}\tilde{e}_{i}+e_{1}\tilde{e}_{1}e_{2}-e_{2}\tilde{e}_{1}e_{1}\right).
\end{align*}
Combining (\ref{t}) with these relations we get
\begin{align*}
\omega_{E}\bullet \mathrm{sym}^{H}\left( \omega_{H}\bullet
{\mathcal T}(\gamma_{0}\alpha_{0})\right) &=
\mathrm{sym}^{H}(hh\tilde{h}) \left(
\beta_{1}\tilde{e}_{1}(e_{1}\wedge e_{2})
+\beta_{2}\tilde{e}_{i}e_{2}e_{i}\right)\\
&+\mathrm{sym}^{H}(hh\tilde{h})\left(\beta_{3}e_{2}(\tilde{e}_{i}\wedge
e_{i})+\beta_{4}e_{i}e_{2}\tilde{e}_{i}\right)\\
&+\beta_{5}\mathrm{sym}^{H}(hh\tilde{h})(e_{2}\tilde{e}_{1}e_{1}-e_{1}\tilde{e}_{1}e_{2})\\
&-\frac{\mathrm{sym}^{H}(hh\tilde{h})}{4n-3}\left(
(4n-1)e_{i}\tilde{e}_{i}e_{2} +\tilde{e}_{i}e_{i}e_{2}\right)\\
&-\frac{2\mathrm{sym}^{H}(hh\tilde{h})}{4n-3}(e_{1}\wedge
e_{2})\tilde{e}_{1} ,
\end{align*}
where the constants $\beta_{i}$ are defined by
\begin{align*}
\beta_{1} =-\frac{2(16n^{2}-4n-1)}{5(4n-3)}\quad \beta_{2}=
-\frac{8n-11}{5(4n-3)}\quad \beta_{3}= -\frac{8n-1}{5(4n-3)}\quad
\beta_{4}=\frac{18n-11}{5(4n-3)}
\end{align*}
and
$$
\beta_{5} = - \frac{2(4n^{2}-11n+11)}{5(4n-3)}.
$$
Skew-symmetrizing
$\omega_{E}\bullet\mathrm{sym}^{H}\left(\omega_{H}\bullet
(\gamma_{0}\alpha_{0})\right)$ in the $E$-variables we  obtain our
claim.
\end{proof}

\begin{cor}\label{lema4}
Both ${\mathcal
T}_{(S^{3}H)\Lambda^{3}_{0}E}(\gamma_{0}\alpha_{0})$ and
${\mathcal T}_{(S^{3}H)E}(\gamma_{0}\alpha_{0})$ are non-zero.
\end{cor}

\begin{proof} Since $\mathrm{pr}_{S^{3}H\Lambda^{3}E}\left(\omega_{E}\bullet
\mathrm{sym}^{H}\left( \omega_{H}\bullet{\mathcal
T}(\gamma_{0}\alpha_{0} )\right) \right)$ is not a multiple of
$\omega_{E}$, ${\mathcal
T}_{(S^{3}H)\Lambda^{3}_{0}E}(\gamma_{0}\alpha_{0})$ is non-zero.
On the other hand, using Lemma \ref{lema3}, it is easy to check
that
\begin{equation}
{\mathcal T}_{(S^{3}H)E}(\gamma_{0}\alpha_{0})
=\frac{4n(n+3)}{4n-3} \mathrm{sym}^{H}(hh\tilde{h})e_{2}.
\end{equation}
\end{proof}

Corollary \ref{lema4} implies Proposition \ref{prop2}. Proposition
\ref{prop1} and Proposition \ref{prop2} conclude the proof of our
main result.

\section{Acknowledgements} I grateful to Paul Gauduchon for many
useful discussions about conformal-Killing forms and to Uwe
Semmelmann for his interest in this work. This work was supported
by Consiliul National al Cercetarii Stiintifice din Invatamantul
Superior, through a CNCSIS grant IDEI "Structuri geometrice pe
varietati diferentiabile", [code 1187/2008].


\begin{thebibliography}{99}


\bibitem{ap} V. Apostolov, D. M. J. Calderbank,
P. Gauduchon: {\it Hamiltonian $2$-forms in K\"{a}hler geometry, I
General Theory}, J. Diff. Geom., vol. 73, no 3 (2006), p. 359-412.


\bibitem{bryant} R. Bryant: {\it Metrics with exceptional
holonomy}, Ann. of Math., vol. 126, no. 2 (1987), p. 525-576.

\bibitem{david} L. David, M. Pontecorvo: {\it A characterization of
quaternionic-projective space by the conformal-Killing equation},
J. London Math. Soc., vol. 80 no. 2 (2009), p. 326-340.


\bibitem{mor}A. Moroianu, U. Semmelmann: {\it Twistor
forms on K\"{a}hler manifolds}, Ann. Sc. Norm. Sup. Pisa, Cl. Sci.
(5) 2 (2003), no. 4, p. 823-845.


\bibitem{sal} S. M. Salamon: {\it Quaternionic-K\"{a}hler manifolds},
Invent. Math. 67, no. 1 (1982), p. 143-171.


\bibitem{sem} U. Semmelmann: {\it Conformal-Killing forms in
Riemannian geometry}, Math. Z. 245 no. 3 (2003), p. 503-527.


\bibitem{swann2} A. Swann: {\it Aspects symplectiques de la
geometrie quaternionique}, C. R. Acad. Paris, t. 308, Serie I,
1989, p. 225-228,

\bibitem{swann3} A. Swann: {\it Some remarks on
quaternion-Hermitian manifolds}, Archivum Math., vol. 33, no. 4
(1997), p. 349-354.



\end{thebibliography}
\end{document}